\documentclass[12pt,french,brochure]{bettyb}

\usepackage{url}

\usepackage{natbib}
  
\makeatletter
  \def\NAT@nmfmt#1{{\scshape\NAT@up#1}}
  \def\bibliography#1{%
    \if@filesw
      \immediate\write\@auxout{\string\bibdata{#1}}%
    \fi
    \expandafter\input{\bbl@main@language bst.tex}%
    \@input@{\jobname.bbl}        
  }
\makeatother

\iffalse
  \usepackage{fontspec}
  \setromanfont[Ligatures=TeX]{TeX Gyre Pagella}
  \usepackage{unicode-math}
\else
  \usepackage[utf8]{inputenc}
  \usepackage{newpxtext,newpxmath}
\fi

\usepackage[T1]{fontenc}
\usepackage[french]{babel}
\let\mathscr \mathcal

\usepackage{hyperref}

\author{Antoine Chambert-Loir}
\address{Université de Paris \\ 
IMJ-PRG \\
F-75013, Paris, France}
\email{antoine.chambert-loir@u-paris.fr}
\title{La conjecture de Mordell: \\ origines, approches, généralisations}

\date{Septembre 2021}
\bbkannee{5\textsuperscript{e} ann\'ee, 2021--2022}
\bbknumero{?}

\begin{abstract}
La conjecture de Mordell prédit qu'une équation diophantienne définissant une courbe projective lisse de genre au moins deux n'a qu'un nombre fini de solutions dans un corps de nombres donné. Le siècle qui s'est écoulé depuis son énoncé, en 1922, a vu plusieurs approches, plusieurs démonstrations, ainsi que de vastes extensions dont la plupart sont encore conjecturales. 
Ce texte, qui reprend l'exposé oral, s'efforce de retracer cette histoire.
\end{abstract}

\alttitle{The Mordell conjecture: origins, approaches, generalizations}
\begin{altabstract}
The Mordell conjecture predicts that a diophantine equation defining
a smooth projective curve of genus at least two has only finity many solutions in a given number field. The century that ran since its statement, in 1922,
gave rise to several approaches, several proofs, and vast extensions most of which are still conjectural.
This text is based on the oral presentation and aims at
recalling this story.
\end{altabstract}

\def\Q{\mathbf Q}
\def\R{\mathbf R}
\def\C{\mathbf C}
\def\Z{\mathbf Z}
\def\P{\mathbf P}
\def\Gal{\operatorname{Gal}}
\def\Hom{\operatorname{Hom}}
\def\End{\operatorname{End}}
\def\tors{{\operatorname{tors}}}

\def\PGL{\operatorname{PGL}}

\def\abs#1{\left\lvert#1\right\rvert}
\def\norm#1{\left\lVert#1\right\rVert}
\begin{document}
\maketitle

\section{Origines}

\subsection{La conjecture de Mordell}

La conjecture dont il sera question dans 
cet exposé naît en 1922,
dans un article \citep{Mordell-1922} consacré aux solutions rationnelles
des équations du troisième ou quatrième degré en deux inconnues. 

Le cas des équations du premier degré est évident,
et l'étude des solutions rationnelles
des équations du second degré relève du théorème de Minkowski 
sur les formes quadratiques (1893).
Plus généralement, une telle équation, de degré quelconque,
définit une courbe algébrique et les articles 
de \cite{HilbertHurwitz-1900} 
et~\cite{Poincare-1901}
semblent être les premiers à avoir
abordé leur arithmétique de ce point de vue géométrique.
Hilbert et Hurwitz se cantonnaient au genre nul,
Poincaré abordait le genre~1 --- les courbes elliptiques ---
en expliquant le procédé de construction de solutions
à partir de sécantes (ou tangentes): si $P,Q$ sont deux points de
la courbe à coordonnées rationnelles, 
la droite qui les joint (la tangente à la courbe
si $P=Q$) coupe la cubique en un troisième point dont
les coordonnées sont également rationnelles.

Néanmoins, au début de son article, Mordell note avec quelque rudesse
à quel point l'état des connaissances est maigre:
\begin{quote}\small
Mathematicians have been familiar
with very few questions for so long a period with so little accomplished in the
way of general results, as that of finding the rational solutions,
or say for shortness, the solutions of indeterminate equations of
genus unity of the forms (...)
\[ y^2=4x^3-g_2x-g_3 \]
and there is no loss of generality in assuming that the coefficients
of all equations in this paper are integers.
\end{quote}

Il y démontre que toutes les solutions rationnelles
peuvent se déduire d'un nombre fini d'entre elles par le procédé
géométrique de cordes et tangentes, 
de manière équivalente, à l'aide des formules d'addition pour les 
fonctions elliptiques.  En termes modernes,
le groupe des points rationnels d'une courbe elliptique
est un groupe abélien de type fini.
Il prouve ainsi un fait que \cite{Poincare-1901}
affirmait banalement, sans démonstration, mais ouvrant ainsi
la voie à l'étude du \emph{rang} de la courbe elliptique,
c'est-à-dire celui du groupe abélien de ses points rationnels.
(Notons que pour Poincaré, le rang était le plus petit nombre
de points de la courbe elliptique à partir on peut en déterminer
tous les autres; c'est un invariant un peu plus compliqué.)

À la fin de son article, 
Mordell  énonce plusieurs problèmes qu'il dit ne pas savoir résoudre:
\begin{itemize}
\item Les problèmes~1, 3, 4, et le 5\textsuperscript e qui les englobe,
forment ce qu'on appelle la \emph{conjecture de Mordell}: 
les équations
diophantiennes définissant
une courbe de genre supérieur ou égal à~$2$ n'ont qu'un nombre
fini de solutions rationnelles.
\item Le problème~2 concerne les 
équations hyperelliptiques de genre au moins un,
de la forme $y^2=f(x)$,
où $f$ est un polynôme à coefficients entiers de degré au moins~$4$;
Mordell demande de prouver qu'une telle équation
n'a qu'un nombre fini de solutions \emph{entières}.
\end{itemize}

\subsection{Le théorème de Weil}
Dans sa thèse, \cite{Weil-1928} étend les résultats de Mordell
à une courbe~$X$ de genre~$g$ quelconque définie sur un corps de nombres~$k$.
Il ne démontra pas la conjecture de Mordell
sur la finitude des points rationnels d'une telle courbe,
et la version qu'il démontre, suggérée par le travail de Poincaré,
concerne les « systèmes rationnels de $g$ points ».

La démonstration de Weil reprend, en la systématisant, celle de Mordell.
En termes modernes, il s'agit de démontrer que le groupe~$J(k)$
des points rationnels de la jacobienne~$J$ de la courbe~$X$
est un groupe abélien de type fini. Weil raisonne uniquement
en termes de diviseurs sur la courbe, qu'il appelle
systèmes virtuels. 

Aux anachronismes près, sa démonstration,
valable pour toute variété abélienne~$A$,
est classiquement séparée en deux étapes:
\begin{itemize}
\item Le théorème faible que le groupe $A(k)/2A(k)$ est fini,
permettant d'écrire un point~$a\in A(k)$ sous la forme
$2a'+b$, où $b$ appartient à un ensemble fini;
\item La théorie des hauteurs exprimant que la
hauteur du point~$a'$ sera, à un terme d'erreur près,
plus petite que celle du point~$a$.
\end{itemize}
La \emph{hauteur} d'un point généralise la notion naïve
de taille d'un nombre rationnel (maximum de son numérateur
et de son dénominateur); comme cette notion naïve,
elle donne lieu à un théorème de finitude:
l'ensemble des points rationnels de hauteur majorée
par toute quantité donnée est fini.

Un argument de descente infinie, pour lequel tant Mordell
que Weil invoquent le nom de Fermat, permet de conclure.

Weil termine son article en posant un certain nombre
de questions arithmétiques sur ce groupe~$J(k)$ qui
reviennent à la détermination de ses facteurs invariants
(rang, torsion) et de son comportement par extension des scalaires.
Il pose aussi la question de les déterminer effectivement.

C'était cependant la conjecture de Mordell que Weil avait en vue,
et même si Hadamard lui avait recommandé 
de ne pas « se contenter d'un demi-résultat », 
le travail de Weil laissait cette conjecture ouverte.

\subsection{Le théorème de Siegel}
C'est à \cite{Siegel-1929} qu'on doit le premier énoncé
de finitude conjecturé par Mordell,
à savoir le cas des solutions entières d'équations
diophantiennes définissant des courbes de genre~$1$ ou plus.
(En fait, Siegel avait déjà obtenu quelques cas  de finitude
dans les années précédentes.)

Le théorème de Siegel vaut aussi
pour une équation diophantienne courbe de genre~$0$ à condition qu'elle ait
au moins trois points à l'infini.
Pour être un peu plus précis, 
une équation $f(x,y)=0$ en deux indéterminées
fournit une courbe affine, à laquelle il manque donc
un certain nombre de points à l'infini. Lorsqu'il s'agit
d'une courbe elliptique dans la présentation $y^2=4x^3-g_2x-g_3$,
seule l'origine de la courbe elliptique est absente.

\cite{Mahler-1934} a généralisé le théorème de Siegel
au cas des $S$-entiers d'un corps de nombres~$k$,
$S$ étant un ensemble fini de places finies de~$k$.
Autrement dit, en langage moderne:

Soit $k$ un corps de nombres, soit $\mathscr O$
un localisé de l'anneau des entiers algébriques par rapport
à un nombre fini d'idéaux maximaux,
soit $X$ est une courbe affine lisse sur~$k$
dont la caractéristique d'Euler est strictement négative;
géométriquement, $X$ est le complémentaire d'un nombre fini, $r$,
de points dans une courbe projective lisse connexe de genre~$g$,
et on suppose l'inégalité $2g-2-r<0$.
Alors pour tout schéma~$\mathscr X$ de type fini sur~$\mathscr O$
qui est un modèle de~$X$,  l'ensemble $\mathscr X(\mathscr O)$
est fini.

Un cas particulier important est l'\emph{équation aux $S$-unités:}
l'équation $x+y=1$ 
n'a qu'un nombre fini de solutions $(x,y)$ dans~$\mathscr O^\times$.

La preuve de Siegel apporte une nouvelle idée,
que d'exploiter les inégalités d'approximations diophantiennes
initiées par Liouville (1844), \cite{Thue-1909}, \cite{Siegel-1921},
\cite{Dyson-1947},
et qui culmineraient avec le théorème de~\cite{Roth-1955}:
\emph{si $\alpha$
est un nombre algébrique irrationnel et~$\rho>1$,
l'ensemble des couples $(p,q)$ tels que $q^\rho \abs{q \alpha  -p} \leq1$
est fini.}
Liouville exigeait l'inégalité $\rho > d-1$, 
où $d$ est le degré de~$\alpha$, Thue traitait
le cas $\rho>d/2$, Siegel et Dyson les cas $\rho > 2\sqrt d$
et $\rho > \sqrt{2d}$.

Lorsque le corps de base est~$k=\Q$,
la première étape de la démonstration de Siegel consiste à
raisonner par l'absurde et à considérer une suite $(a_n)$
de points de points entiers de~$X$
que, par compacité, on peut supposer converger, pour la topologie de~$\R$,
vers un des points à l'infini de la courbe,
en particulier vers un point algébrique. L'intégralité des points~$a_n$
fournit une inégalité dans un sens opposé à ce qu'interdisent
les résultats d'approximation diophantienne,
mais toutefois trop faible pour les contredire.

Dans une seconde étape, Siegel utilise les propriétés de division,
via des revêtements non ramifiés appropriés,
qui permettent d'améliorer l'inégalité d'approximation,
jusqu'à contredire les théorèmes de Thue ou, plus précisément,
le raffinement qu'il avait démontré.

\section{Approches}

\subsection{La méthode de Chabauty}

La méthode de \cite{Chabauty-1941} est peut-être l'approche
la plus simple de la conjecture de Mordell, mais elle fonctionne
uniquement sous l'hypothèse restrictive que
le rang de Mordell–Weil de la courbe~$X$ considérée
soit strictement inférieur à son genre.

Dans ce cas, \cite{Chabauty-1941} fixe un nombre premier auxiliaire~$p$
et considère l'adhérence~$\Gamma$ dans~$J(\Q_p)$, pour la topologie $p$-adique,
du groupe $J(\Q)$. Si $r$ est le rang de~$J(\Q)$
et $g$ le genre de~$X$, donc la dimension de~$J$,
cette adhérence est un sous-groupe de Lie $p$-adique
de dimension~$r<g$.

Il existe donc une forme différentielle invariante~$\omega$ sur~$J$,
non nulle, dont la restriction à~$J(\Q_p)$, est nulle.
On peut définir, si $a$ et~$b$ sont des points de~$J(\Q_p)$
suffisamment proches, l'intégrale $\int_a^b \omega$; 
il suffit d'intégrer formellement dans une carte où tout cela fait sens.
Si les points~$a$ et~$b$ appartiennent à~$\Gamma$, cette intégrale sera nulle.

D'autre part, si $a$ et~$b$ appartiennent à~$X(\Q_p)$,
ces intégrales sont celles de la forme différentielle $\omega|_X$.
Comme le plongement d'Abel–Jacobi de~$X$ dans~$J$ induit une bijection
de~$\Omega^1_J$ sur~$\Omega^1_X$, la forme $\omega|_X$ n'est pas nulle.
Par suite, le point~$a$ étant  fixé dans~$X(\Q)$,
l'intégrale $\int_a^x\omega$,
qui est une fonction analytique locale sur~$X(\Q_p)$,
n'a qu'un nombre fini de zéros proches de~$a$.
Cela entraîne que $X(\Q)$ n'a qu'un nombre fini de points proches de~$a$.

Par compacité de~$X(\Q_p)$, il en résulte que $X(\Q)$ est fini.

Cette méthode a été précisée par~\cite{Coleman-1985}
qui propose une majoration effective, et souvent efficace,
du nombre de points rationnels de~$X$.
 
Elle a connu de nombreux développements dans les dernières
années. 
La version non linéaire de~\cite{Kim-2005,Kim-2006}
montre un lien entre des conjectures « motiviques »
comme la conjecture de Bloch–Kato et la conjecture de Mordell.
Si elle ne permet pas de démontrer la conjecture
de Mordell en toute généralité, elle donne lieu 
à une nouvelle démonstration du théorème de Siegel.
Une série d'articles remarquables
\citep{BalakrishnanDogra-2018,BalakrishnanDogra-2021,BalakrishnanBesserBianchiEtAl-2021} parvient à mettre en œuvre une partie de la méthode de Kim:
c'est ce qu'on appelle la méthode « quadratique » de Chabauty.
Elle a notamment permis à \cite{BalakrishnanDograMullerEtAl-2019}
de déterminer les points rationnels
de la courbe modulaire $X_{\mathrm {split}}(13)$,
finissant ainsi la preuve par \cite{BiluParent-2011,BiluParentRebolledo-2013}
du cas « sous-groupe de Cartan déployé » d'une conjecture
de Serre sur l'image de la représentation galoisienne associée à une courbe
elliptique sur~$\Q$.

\subsection{Le cas des corps de fonctions}

La géométrie diophantienne « sur un corps de fonctions »
cherche à étudier des questions analogues  à des problèmes
arithmétiques connus en remplaçant
le corps des nombres rationnelles  par le corps
des fonctions d'une courbe. 
Ce glissement est motivé par la proximité algébrique
des deux types de corps --- 
ils sont par exemple tous deux le corps des fractions
d'anneaux de Dedekind --- et par la possibilité de mettre
en œuvre les outils de la géométrie algébrique ou analytique.

Soit $B$ une courbe projective lisse sur un corps~$k$
et soit $X$ une courbe projective lisse sur le corps~$K=k(B)$
des fonctions rationnelles sur~$B$.
À une réserve près que nous discuterons bientôt,
l'analogue de la conjecture de Mordell
demande par exemple de prouver que l'ensemble $X(K)$ est fini.
Adoptant un point de vue complètement géométrique, on peut considérer une
surface lisse~$S$ sur~$k$, fibrée en courbes sur~$B$,
dont la fibre générique est précisément la courbe~$X$.
Les points rationnels de~$X(K)$ correspondent alors
à des \emph{sections} $B\to S$ de la projection canonique.

Le problème contient ainsi, comme cas particulier,
une question géométrique très simple. Supposons en effet
que la surface~$S$
soit une famille de courbes essentiellement constante, 
birationnelle à un produit $B\times C$ de la courbe~$B$
par une courbe~$C$, dont le genre est au moins~$2$.
Les sections $B\to S$ correspondent alors aux morphismes de~$B$
dans~$C$. La finitude de l'ensemble de ceux \emph{qui ne sont
pas constants} est le théorème de de Franchis,
voir \citep[chap.~21, th.~8.27]{ArbarelloCornalbaGriffiths-2011}.
 Là est la réserve nécessaire à l'énoncé de la conjecture
de Mordell sur les corps de fonctions: si l'ensemble des points
rationnels est infini, 
c'est que la courbe initiale donne lieu à une famille
essentiellement constante et sauf un nombre fini d'entre eux,
ces points rationnels
correspondent à des sections constantes.

On doit à \cite{Manin-1963} et \cite{Grauert-1965}
des démonstrations de la conjecture de Mordell sur les corps
de fonctions, en caractéristique zéro.
Le cas de la caractéristique positive a été étudié par \cite{Samuel-1966}.
Ces preuves considèrent la différentielle d'une section~$f$ de~$B$
dans~$S$, c'est-à-dire l'application $f^*\Omega^1_S\to\Omega^1_B$
et exploitent que le genre de~$X$ est au moins~$2$
via la positivité du fibré relatif $\Omega_{S/B}$.

Peu après, \cite{Parsin-1968} a proposé une troisième
démonstration, en caractéristique zéro, mais de nature fondamentalement
différente.

\subsection{La conjecture de Shafarevich et la construction de Kodaira–Parshin}

Si l'on écrit une courbe elliptique sous la forme
d'une équation $y^2=4x^3-g_2x-g_3$, 
son discriminant est défini par $\Delta=g_2^3-27 g_3^2$.
Lorsque $g_2$ et~$g_3$ appartiennent à un anneau et que $\Delta$
y est inversible, la courbe elliptique a alors bonne réduction partout.
Ainsi, prenant $g_2$ et~$g_3$ dans l'anneau des entiers d'un corps de nombres,
la courbe aura bonne réunion hors des idéaux premiers 
qui contiennent~$\Delta$.

\cite{Shafarevich-1963} avait déduit du théorème de Siegel 
un énoncé de finitude : un corps de nombres~$k$ et un ensemble
fini de places de~$k$ étant donnés, l'ensemble
des courbes elliptiques sur~$k$ qui ont bonne réduction
hors des places de~$S$ est fini.

Il conjectura alors le résultat analogue pour les variétés
abéliennes de dimension arbitraire,
ainsi que pour les courbes projectives de genre arbitraire.
Par la construction de la jacobienne d'une courbe,
l'étude de sa mauvaise réduction, et
 le théorème de Torelli, le résultat pour les variétés
abéliennes entraîne le cas des courbes.
Il conjecture également l'inexistence
de courbes projectives lisses et de variétés abéliennes
sur~$\Q$ qui aient bonne réduction partout, 
une conjecture qui sera démontrée par~\cite{Fontaine-1985}.

\cite{Parsin-1968} a montré comment déduire la conjecture
de Mordell de celle de Shafarevich. Sa méthode
consiste à construire,
pour toute courbe~$B$ de genre~$\geq 2$, une famille
non isotriviale de courbes $S\to B$ paramétrée par~$B$;
en fait, si $b\in B$, la courbe~$S_b$ est construite
comme un revêtement étale de~$B$ ramifié en~$b$ uniquement.
Par suite, pour tout~$b\in B$, il n'y a qu'un nombre fini
de $c\in B$ tels que les courbes~$S_b$ et~$S_c$ soient isomorphes.
Une construction similaire avait également proposée par \cite{Kodaira-1967}.

Lorsque $b$ varie dans~$B(k)$, cela fournit une famille de courbes sur~$k$,
chacune n'étant isomorphe qu'à un nombre fini de membres de la famille.
en faisant cette construction de façon arithmétique,
on constate que les courbes~$S_b$ n'ont mauvaise réduction
qu'en un ensemble fini de places de~$k$. D'après 
la conjecture de Shafarevich, l'ensemble~$B(k)$ est fini.

Pour être précis, la construction requiert une extension finie~$K$
du corps de base ainsi que de remplacer~$B$ par un revêtement
étale~$B_1$. D'après le théorème de Chevalley–Weil, lui-même
une conséquence du théorème de finitude de Hermite–Minkowski,
il existe une extension finie~$k_1$ de~$K$ telle que
$B_1(k_1)$ contienne l'image réciproque de~$B(k)$ dans~$B_1$.
La discussion précédente entraîne que $B_1(k_1)$ est fini,
donc $B(k)$ aussi.

Par ailleurs, \cite{Parsin-1968} et~\cite{Arakelov-1971}
en caractéristique zéro, et \cite{Szpiro-1979} en caractéristique
positive, démontrent l'analogue
de la conjecture de Shafarevich 
sur les corps de fonctions de caractéristique zéro, 
ce qui fournit ainsi une troisième démonstration de la conjecture
de Mordell dans ce cadre.

\subsection{La démonstration de Faltings}

C'est à \cite{Faltings-1983a} qu'on doit la première démonstration
de la conjecture de Mordell. Il suit l'idée
de Parshin et Arakelov en ce qu'il démontre d'abord
la conjecture de Shafarevich, en la réduisant à une autre
conjecture, proposée par \cite{Tate-1966}, sur la représentation galoisienne
$\ell$-adique associée à une variété abélienne.

Soit $A$ une variété abélienne sur un corps~$k$, soit $g$ sa dimension,
et soit $\ell$ un nombre premier distinct de la caractéristique de~$k$.
Pour tout entier~$n\geq 1$, le groupe~$A_{\ell^n}$ 
des points~$a\in A(\overline k)$
tels que $\ell^n a=0$ est isomorphe à~$(\Z/\ell^n\Z)^{2g}$,
et est muni d'une action du groupe $G=\Gal(\overline k/k)$ 
Le module de Tate $\ell$-adique est la limite $T_\ell(A)\simeq \Z_\ell^{2g}$
des groupes finis~$A_{\ell^n}$; l'action de~$G$ sur les
groupes finis~$A_{\ell^n}$ induit une action continue de~$G$
sur $T_\ell(A)$. 
On peut d'ailleurs se contenter d'étudier la représentation
galoisienne sur l'espace vectoriel $V_\ell(A)=\Q_\ell\otimes_{\Z_\ell} T_\ell(A)$.

La construction de $V_\ell(A)$ est fonctorielle.
Si $f\colon A\to B$ est un homomorphisme de variétés abéliennes
sur~$k$, elle donne lieu à un morphisme $f_*$, $\Q_\ell$-linéaire,
$G$-équivariant, de $V_\ell(A)$ dans~$V_\ell(B)$.
La \emph{conjecture de Tate} est que cette fonctorialité donne lieu 
à une {bijection}
\[ \Q_\ell \otimes_\Z \Hom_k(A,B) \to \Hom_G(V_\ell(A),V_\ell(B)). \]
Comme la $\Q$-algèbre $\Q\otimes_\Q\End_k(A)$
déduite de celle des endomorphismes d'une variété abélienne
est semi-simple, une conséquence de cette conjecture de Tate
est la semi-simplicité de la représentation galoisienne~$ V_\ell(A)$.

Le théorème principal de~\cite{Tate-1966} est la démonstration
de cette conjecture lorsque $k$ est un corps fini.
Un argument essentiel est la finitude de l'ensemble
de classes d'isomorphismes de variétés abéliennes sur~$k$
munies d'une polarisation de degré donné et d'une isogénie vers~$A$
de degré une puissance de~$\ell$. Comme le fait remarquer~\cite{Faltings-1986},
cette démonstration déduit la conjecture de Tate
de la conjecture de Shafarevich, alors que \cite{Faltings-1983a}
déduit la conjecture de Shafarevich de la conjecture de Tate.

Il n'est pas possible ici de résumer fidèlement cette démonstration de Faltings.
Disons seulement que trois théorèmes de finitude interviennent.
Le premier repose sur la notion de hauteur que Faltings introduit sur l'espace
des modules des variétés abéliennes et sur la finitude
de l'ensemble des classes d'isomorphie de variétés abéliennes 
définies sur un corps de nombres donné,
munie d'une polarisation de degré donné,
et de hauteur bornée. 
En étudiant le comportement de cette hauteur par isogénie,
Faltings en déduit la finitude de l'ensemble
des classes d'isomorphie de variétés abéliennes sur un corps de nombres
qui isogènes à  une variété abélienne donnée puis,
par la méthode de Tate, la conjecture de Tate elle-même.

Le troisième théorème de finitude que intervient pour démontrer
la conjecture  de Shafarevich est  celui de Hermite–Minkowski,
à savoir la finitude de l'ensemble des corps de nombres de degré
donné et qui sont non ramifiés en dehors d'un ensemble fini
de nombres premiers donné. 

Soit alors $k$ un corps de nombres, soit $S$ un ensemble fini de places de~$k$
et soit~$A$ une variété abélienne de dimension~$g$ définie sur~$k$
et ayant bonne réduction hors des places de~$S$. Choisissons un nombre
premier~$\ell$ et adjoignons à~$S$, si nécessaire, toutes les places
divisant~$\ell$. Du théorème de Hermite–Minkowski,
du critère de Néron–Ogg–Shafarevich qui assure que
la représentation~$V_\ell(A)$ est non ramifiée hors des places de~$S$,
et du théorème de densité de Čebotarev, on déduit
--- voir \cite[th.~3.1]{Deligne-1985} ---
que la représentation galoisienne~$V_\ell(A)$ est
déterminée par les traces d'un nombre fini d'endomorphismes de Frobenius
agissant sur~$V_\ell(A)$.
Compte tenu du théorème de Weil qui les majore par $2g\sqrt q$
(si $q$ est le cardinal du corps résiduel de ces endomorphismes
de Frobenius), il n'y a qu'un nombre fini de représentations
galoisiennes possibles pour $V_\ell(A)$.

\subsection{La démonstration de Vojta}

L'article de \cite{Vojta-1991} propose une nouvelle preuve, totalement
différente, de la conjecture de Mordell, reposant sur une inégalité de
hauteurs.
 
Pour évoquer cette inégalité, il faut rappeler
que la hauteur des points d'une variété abélienne~$A$ définie
sur un corps de nombres~$k$, hauteur dont Mordell et Weil avaient
fait un usage fondamental, peut être promue en une forme quadratique
définie positive sur le $\R$-espace vectoriel
$\R\otimes_\Z A(k)$, espace vectoriel qui est de dimension finie
en raison du théorème de Mordell–Weil.
Par la fonctorialité approchée des hauteurs, cette hauteur
est une « forme quadratique approchée », mais
Néron et Tate avaient construit une forme quadratique canonique
qui en diffère par une fonction bornée. Cette forme
quadratique, qu'on appelle maintenant \emph{hauteur de Néron–Tate},
est a priori positive sur l'espace vectoriel réel~$\R\otimes_\Z A(k)$,
mais on déduit du théorème de finitude de Northcott qu'elle
est effectivement définie positive. De plus, l'homomorphisme canonique
de~$A(k)$ dans $\R\otimes_\Z A(k)$ a pour noyau  le sous-groupe~$A_\tors$,
et le théorème de Mordell–Weil (ou, plus directement, celui de Northcott)
entraîne également que ce sous-groupe est fini.

Soit $X$ une courbe projective lisse de genre~$g\geq 2$
définie sur un corps de nombres, considérée comme plongée
dans sa jacobienne~$J$ --- au moyen d'un point rationnel arbitraire.
Posons $V=\R\otimes_\Z J(k)$ et soit $\langle\cdot,\cdot\rangle$
et $\norm\cdot$ 
le produit scalaire et la norme de Néron–Tate sur~$V$.
Le théorème de \cite{Vojta-1991} est une inégalité pour l'angle
$\angle(x,y)$ que forment deux points~$x$ et~$y$
de~$X(k)$ dans l'espace vectoriel~$V$:
il existe des nombres réels~$c$ et~$c'$
tels que si $x,y\in X(k)$ satisfont $\norm x\geq c$ et
$ \norm y\geq c' \norm x$, alors $\cos(\angle(x,y))\geq 3/4$.

La conjecture de Mordell en résulte par un argument géométrique très simple.
On commence par recouvrir l'espace vectoriel~$V$ par un nombre fini de cônes
d'angle strictement inférieur à $\arccos(3/4)$ et l'on prouve
que chacun de ces cônes~$C$ ne contient qu'un nombre fini de points de~$X(k)$.
Si tous les points~$x \in X(k)$ de~$C$ vérifient
$\norm x\leq c$, alors~$C$ ne contient qu'un nombre fini de points de~$X(k)$.
Supposons que~$C$  contienne  un point~$x$ tel que $\norm x\geq c$;
d'après l'inégalité de Vojta,
tous les autres points~$y$ de~$X(k)$ contenus dans ce cône
satisfont $\norm y\leq c'\norm x$; il s'ensuit que $C$ ne contient
qu'un nombre fini de points de~$X(k)$.

Signalons que \cite{Mumford-1965} avait établi
une inégalité de hauteurs similaire: si $a>1/2g$, 
il existe un nombre réel~$c$ tel que  pour tous $x,y\in X(k)$
tels que $x\neq y$, on ait
$\langle x,y\rangle \leq a( \norm x^2+\norm y^2)+c$,
soit encore $\norm{x-y}^2\geq (1-2a) (\norm x^2+\norm y^2)-2c$.
Chez Mumford, cette inégalité est utilisée comme un principe d'espacement,
il en déduit que parmi les $\approx T^r$ points~$x$ de~$J(k)$
tels que $\norm x\leq T$, au plus $\approx \log(T)$ appartiennent à~$X(k)$.

Dans le titre de l'article de~\cite{Vojta-1991},
l'allusion au théorème de Siegel
montre qu'il s'agit de nouveau de techniques approximation diophantienne.
Dans cet article, ces méthodes sont exprimées
dans le langage de la géométrie d'Arakelov,
telle que développée par~\cite{Arakelov-1974b}, \cite{Faltings-1986},
\cite{GilletSoule-1990}.
\cite{Bombieri-1990} en a aussitôt donné une présentation 
dans le langage de la géométrie diophantienne classique.

Au dialecte utilisé près,
toutes les démonstrations d'approximation diophantienne
depuis Liouville
suivent un schéma similaire. Dans une première
phase, on construit un polynôme auxiliaire
en les données, non nul, qui s'annule sur un ensemble prescrit avec
des multiplicités prescrites, dont les coefficients
sont de taille contrôlée. C'est ici qu'interviennent
le \emph{lemme de Siegel} ou, chez Vojta, 
le théorème d'« amplitude arithmétique » de \cite{GilletSoule-1988}.
Il s'agit alors d'étudier une dérivée convenable de ce polynôme auxiliaire,
a priori non nulle, et d'utiliser le « théorème fondamental de l'arithmétique » selon lequel un entier positif non nul   est au moins égal à~$1$.
Garantir cette non-annulation est évident dans l'inégalité de Liouville,
mais l'est beaucoup moins pour les développements ultérieurs.
Vojta utilise un énoncé géométrique \citep{Vojta-1989}
qu'il avait
démontré pour donner une démonstration \citep{Vojta-1989a}
de la conjecture de Mordell sur les corps de fonctions.
Bombieri utilise un énoncé arithmétique plus élémentaire,
le « lemme de Roth » de~\cite{Roth-1955}.

\subsection{La démonstration de Lawrence et Venkatesh}

C'est la dernière en date \citep{LawrenceVenkatesh-2020},
nous n'en dirons rien ici en renvoyant à l'exposé
de Marco Maculan au Séminaire Bourbaki.

\section{Généralisations}

\subsection{Dimension supérieure}

On doit à Lang d'avoir posé la question 
d'une généralisation de la conjecture de Mordell pour les variétés
de dimension arbitraire, d'abord pour des variétés abéliennes
dans~\citep{Lang-1960a}, et surtout en général dans~\citep{Lang-1986a}.

Il y a en fait 
deux façons d'étendre la conjecture de Mordell en dimension supérieure:
on peut demander que l'ensemble des points rationnels soit fini,
ou bien simplement qu'il ne soit pas dense pour la topologie de Zarisi.

De même, si $X$ est une variété algébrique sur un corps~$k$,
il y a plusieurs façons d'étendre l'hypothèse « $g\geq 2$ » 
dans la conjecture de Mordell:
\begin{itemize}
\item Si $X$ est lisse, que $\omega_X$ soit ample (Lang dit que $X$ est \emph{canonique});
\item Quitte à considérer une désingularisation  de~$X$ est lisse, 
que $\omega_X$ (et ses puissances) définisse une application
birationnel vers un espace projectif ($X$ est de \emph{type général};
Lang dit que $X$ est \emph{pseudo-canonique});
\item Qu'il n'existe aucun morphisme non constant de~$\P_1$ ou d'une
variété abélienne vers~$X$;
\item Si $k$ est un sous-corps de~$\C$, qu'il n'existe aucune application
holomorphe de~$\C$ dans~$X(\C)$ (c'est-à-dire que $B$ soit \emph{hyperbolique au
sens de Brody});
\item Toujours si $k$ est un sous-corps de~$\C$,
que  $X(\C)$ soit \emph{hyperbolique au sens de Kobayashi.}
\end{itemize}

L'article de~\cite{Lang-1986a} discute ces variantes,
et leurs relations possibles avec les points rationnels.
On peut espérer les énoncés suivants, pour une variété
projective lisse sur un corps de nombres~$k$:
\begin{itemize}
\item Si $X$ est de type général, l'ensemble $X(k)$ n'est
pas dense pour la topologie de Zariski,
un énoncé également proposé par Bombieri;
\item Si $X(\C)$ est hyperbolique, l'ensemble $X(k)$ est fini.
\end{itemize}
Ces énoncés sont encore complètement ouverts,
et en suggèrent de nombreux autres. Par exemple,
si $X$ est de type général, quelle est la réunion
des sous-variétés de dimension strictement positive
de~$X$ qui ne sont pas de type général?

On doit à \cite{Campana-2004} d'avoir proposé une 
description conjecturale complète de ce que
pourrait être l'adhérence de~$X(k)$ pour la topologie de Zariski
lorsque $X$ est une variété projective~$X$ arbitraire
sur un corps de nombres~$k$. 

\cite{Vojta-1987} a proposé un dictionnaire entre 
géométrie arithmétique et théorie de Nevanlinna qui, en particulier,
fait correspondre le théorème de \cite{Roth-1955}
avec le deuxième théorème principal de Nevanlinna.
La théorie de Nevanlinna en dimension supérieure de Griffiths 
suggère une conjecture remarquable d'approximation diophantienne
en toute dimension dont Vojta a démontré qu'elle entraîne 
non seulement la conjecture de Bombieri–Lang que les points
rationnels d'une variété de type général ne sont pas
denses pour la topologie de Zariski, mais également
des énoncés analogues pour les points entiers.

Enfin, ces questions peuvent être posées, non seulement
sur un corps de nombres, mais sur tout corps de type fini,
voire sur tout corps de type fini sur un corps algébriquement clos.

\subsection{La conjecture de Zilber–Pink}

Il y a cependant un cas des conjectures de Bombieri–Lang
où l'on sait quelque chose,
encore une fois grâce à \cite{Faltings-1991,Faltings-1994}: lorsque
$X$ est une sous-variété d'une variété abélienne~$A$ définie
sur un corps de nombres~$k$,
l'adhérence de~$X(k)$ pour la topologie de Zariski
est une réunion finie de translatées de sous-variétés abéliennes
de~$A$. Il faut mettre cet énoncé en relation
avec le théorème de~\cite{Kawamata-1980}
qui affirme qu'une sous-variété d'une variété abélienne
est de type général si et seulement si ce n'est pas un translaté
d'une sous-variété abélienne.

Plus généralement, \cite{Hindry-1988} remplace~$X(k)$
par l'intersection de~$X(\overline k)$
avec un sous-groupe~$\Gamma$ de rang fini de~$A(\overline k)$,
c'est-à-dire tel que $\Q\otimes_\Z\Gamma$ soit un $\Q$-espace
vectoriel de dimension finie.

Ces théorèmes ont été généralisés par \cite{McQuillan-1995}
au cas des sous-variétés de variétés semi-abéliennes.

Par ailleurs, \cite{Remond-2000a} a démontré l'analogue de l'inégalité
de Vojta dans ce contexte.

Pour conclure, mentionnons que \cite{Zilber-2002} et \cite{Pink-2005a,Pink-2005}, de façon indépendante, ont émis des conjectures
sur les intersections atypiques d'une sous-variété d'une variété
semi-abélienne avec des sous-groupes algébriques propres.
\emph{Atypique} signifie ici que l'intersection est de dimension 
strictement supérieure à ce qu'un décompte de paramètre permet d'escompter.
En fait, l'énoncé proposé par Pink prend en compte les variétés
de Shimura mixtes --- il impliquerait également la conjecture d'André–Oort.
Ce sujet a connu une intense activité dans ces dernières
années, je renvoie à la synthèse de \cite{HabeggerRemondScanlonEtAl-2017}.

\subsection{Uniformité}

Rappelons que la méthode de Chabauty–Coleman, lorsqu'elle s'applique,
fournit une majoration explicite du nombre de points rationnels
d'une courbe de genre au moins~$2$.

Dès la preuve par \cite{Faltings-1983a} de la conjecture de Mordell,
il était apparu possible
d'en déduire des informations plus précises.
Si $X$ est une courbe projective lisse de genre~$g$
définie sur un corps de nombres~$k$,
l'exposé~XI de~\cite{Szpiro-1985} propose ainsi une majoration
du cardinal de~$X(k)$ 
qui fait intervenir un certain nombre d'expressions,
dont le rang du groupe de Mordell–Weil $J(k)$ de la jacobienne de~$X$.

La démonstration de~\cite{Vojta-1991} donne également lieu
à de telles majorations. Leur intérêt est d'être relativement
uniforme (à l'utilisation du rang de Mordell–Weil près) 
lorsque la courbe~$X$ varie dans une famille. 
Cela était par exemple indiqué par un théorème de~\cite{deDiego-1997}
qui, cependant, ne décomptait  que les points de hauteur assez grande.

Dans le contexte de la conjecture de Lang
sur les sous-variétés de variétés abéliennes,
\cite{Remond-2000b} a majoré de façon explicite le nombre
de translatés de sous-variétés de variétés abéliennes et
les degrés de ces sous-variétés qui interviennent par
une expression ne dépendant que de la dimension du degré de la sous-variété,
la dimension de la variété abélienne ambiante, et le rang du groupe~$\Gamma$.

Tout récemment, \cite{DimitrovGaoHabegger-2021} ont démontré
une majoration du cardinal de~$X(k)$ par une expression
de la forme $c^{1+\rho}$, où $c$ est une constante ne dépendant
que du genre de~$X$ et du degré du corps de nombres~$k$,
et $\rho$ est le rang de~$J(k)$.

Peut-on aller plus loin ? C'est bien possible. Même si elle
semble largement inaccessible, une piste consisterait à majorer
uniformément le rang de~$J(k)$; il n'est en fait même pas
clair que ce soit vrai.
Toutefois, 
sous l'hypothèse que la conjecture de Bombieri–Lang est vérifiée,
\cite{CaporasoHarrisMazur-1997} et \cite{Pacelli-1997}
ont démontré
que sur un corps de nombres~$k$ de degré donné,
le nombre de points rationnels d'une courbe de genre~$g$
est uniformément borné ! Leur démonstration repose
sur un énoncé géométrique inconditionnel: l'existence,
pour une famille $f\colon X\to B$ de courbes, propre, génériquement lisse,
d'un entier~$n$ tel que le produit symétrique fibré $S^n_B X$ admette
une application rationnelle dominante vers une variété~$W$ de type général.
Appliquée à~$W$, la conjecture de Bombieri–Lang contraint
les points rationnels de~$S^n_B X$ (qui sont, en famille,
les diviseurs effectifs de degré~$n$ d'une même courbe $X_b$,
pour $b\in B(k)$) à être contenus dans une sous-variété stricte.
En partant d'une famille (uni)verselle de courbes de genre~$g$,
un argument de récurrence leur permet de conclure.

\subsection{Variantes de la conjecture de Shafarevich}

Il s'agit de démontrer, si possible,
une finitude de l'ensemble des classes d'isomorphie
de variétés d'un type géométrique donné, définies
sur un corps de nombres, et ayant bonne réduction
en dehors d'un ensemble fini de places donné.

Du théorème de Faltings pour les variétés abéliennes,
\cite{Andre-1996} a par exemple déduit le cas des surfaces~K3
(plus généralement, de certaines classes de variétés hyperkählériennes),
tandis que \cite{JavanpeykarLoughran-2017} traitaient
le cas des hypersurfaces de niveau de Hodge~$\leq 1$,
ou de certaines surfaces de type général.

L'approche de \cite{LawrenceVenkatesh-2020} permet
d'établir de nouveaux cas.
Par exemple, si $n$ est assez grand et
$d$ est assez grand (dépendant de~$n$), alors pour tout entier~$N\geq 1$,
dans l'espace projectif des hypersurfaces
de degré~$d$ de~$\P_n$, celles qui ont bonne réduction
en toute place ne divisant pas~$N$ ne sont pas denses
pour la topologie de Zariski.
La conjecture de Lang–Vojta entraîne que cet ensemble
est fini modulo l'action de $\PGL(n+1,\Z[1/N])$.

\subsection{Effectivité}

Pour terminer sur une note spéculative, il faut mentionner
la question de l'effectivité. Dans le cadre
de la conjecture de Mordell, cela voudrait dire un moyen
pratique de déterminer l'ensemble des points rationnels
d'une courbe de genre~$\geq2$ définie sur un corps de nombres,
dans celui du théorème de Siegel,
l'ensemble des points entiers de la courbe en question,
dans celui du théorème de Mordell–Weil, 
un ensemble générateur du groupe
des points rationnels d'une variété abélienne définie sur un corps 
de nombres; dans celui de la conjecture de Shafarevich,
l'ensemble des variétés abéliennes de dimension donnée
définies sur un corps de nombres donné ayant bonne réduction
en dehors d'un ensemble fini donné de places.

Bien sûr, on peut tenter d'énumérer les solutions
mais,
même en présence de majorations du nombre de ces solutions,
comment savoir que la recherche est achevée?

Le point de vue des hauteurs offre une formulation possible:
il suffirait de démontrer une majoration explicite (ou explicitable)
pour la hauteur des points en question, pour la hauteur de Faltings
des variétés abéliennes en question.

Même s'il  faut mentionner que la théorie
des formes linéaires  en logarithmes de \cite{Baker-1975}
permet de traiter le cas des points entiers
des courbes hyperelliptiques,
les résultats sont hélas très partiels.

En ce qui concerne le groupe de Mordell–Weil,
\cite{Tate-1974} propose un algorithme, conditionnel à
la finitude du groupe de Tate–Shafarevich, 
finirait par le déterminer explicitement.

\cite{Moret-Bailly-1990} introduit une version
effective relativement explicite de la conjecture de Mordell:
si $X$ est une courbe de genre~$g\geq 2$ définie sur un corps de nombres~$k$,
la hauteur de tout point algébrique~$P$ de~$X$
devrait être majorée par une expression de la forme
\[ A([k(P):k]) \log(\operatorname{disc}(k(P)) + B([k(P):k]), \]
où $A$ est une fonction dépendant de la définition choisie
pour la hauteur sur~$X$.

Cet article explicite également
les interactions entre cette conjecture et d'autres
énoncés de géométrie diophantienne,
la conjecture de Szpiro sur le discriminant minimal des courbes
elliptiques et la conjecture~ABC de Masser et Oesterlé,
voir aussi \cite{Oesterle-1988}.
Indépendamment, \cite{Elkies-1991a} avait montré que cette
dernière conjecture implique la conjecture de Mordell.

\section*{Indications bibliographiques}

Il y a trois magnifiques ouvrages d'introduction
à la géométrie diophantienne, ceux de
\cite{Serre-1997a}, \cite{HindrySilverman-2000}
et \cite{BombieriGubler-2006}. Le premier, antérieur
à la démonstration de Faltings, couvre néanmoins
le théorème de Mordell–Weil, le théorème de Siegel
et les théorèmes de Chabauty et Mumford;
les deux autres présentent en outre la démonstration
du théorème de Roth et 
celle de la conjecture de Mordell par Bombieri.

Le survol de \cite{Lang-1991} aborde la première démonstration
de la conjecture de Mordell par Faltings,
ainsi que les conjectures en dimension supérieure.
Une analyse détaillée de cette démonstration fait l'objet
des comptes rendus de séminaires de \cite{Szpiro-1985}
et \cite{FaltingsWustholz-1992}.
Le séminaire de~\cite{Szpiro-1990} est guidé par l'étude
d'une version effective de la conjecture de Mordell.

\bibliographystyle{mynat}
\itemsep 0pt plus 1pt
\bibliography{aclab,zot-mordell}
\end{document}